\documentclass{amsart}
\newcommand\Quad{\HS{0.5}}

\usepackage{amssymb, amsmath}
\usepackage{graphics}
\usepackage{pst-all}
\usepackage{color}
\usepackage[english,francais]{babel}

\input xypic
\xyoption{all}
\xyoption{arc}

\newtheorem{thrm}{Theorem}[section]
\newtheorem{coro}[thrm]{Corollary}
\newtheorem{defi}[thrm]{Definition}
\newtheorem{lemm}[thrm]{Lemma}
\newtheorem{Eprop}[thrm]{Proposition}

\newenvironment{Proof}
{\noindent\emph{Proof.}}
{\null\hfill$\square$\medskip}

\newenvironment{sketch}
{\noindent\emph{Proof (sketch).}}
{\null\hfill$\square$\medskip}

\psset{unit=1mm}
\psset{fillcolor=white}
\psset{dotsep=1.5pt}
\psset{dash=1.5pt 1.5pt}
\psset{linewidth=0.8pt}
\psset{arrowsize=3.5pt}
\psset{doublesep=1pt}
\psset{coilwidth=1.2mm}
\psset{coilheight=2}
\psset{coilaspect=20}
\psset{coilarm=2}
\newpsstyle{double}{linewidth=0.5pt, doubleline=true}
\newpsstyle{etc}{linestyle=dotted}
\newpsstyle{exist}{linestyle=dashed}
\newpsstyle{thin}{linewidth=0.5pt}
\newpsstyle{thinexist}{linewidth=0.5pt, linestyle=dashed}
\def\AArrow(#1,#2){\ncline[nodesep=0.5mm, linewidth=0.8pt, border=1.2pt]{->}{#1}{#2}}
\def\BArrow(#1,#2){\ncline[linecolor=blue, linewidth=1.2pt, linestyle=dotted, nodesep=0.5mm, border=1.2pt]{->}{#1}{#2}}
\def\CArrow(#1,#2){\ncline[linecolor=red, nodesep=0.5mm, linewidth=0.8pt, border=1.2pt]{->}{#1}{#2}}

\newcounter{ITEM}
\newcommand\ITEM[1]{\setcounter{ITEM}{#1}\leavevmode\hbox{\rm(\roman{ITEM})}}
\newcommand\card{\mathtt{\#}}
\newcommand\cone{\operatorname{cone}}
\newcommand\dive{\preccurlyeq\nobreak}
\newcommand\domby{\preccurlyeq\nobreak}
\newcommand\ff{f}
\newcommand\Gar{F}
\newcommand\Garext{E}
\renewcommand\ge{\geqslant}
\renewcommand\gg{g}
\newcommand\hh{h}
\newcommand\HS[1]{\hspace{#1ex}}
\newcommand\ie{\textit{i.e.}}
\newcommand\ii{i}
\newcommand\II{I}

\newcommand\kk{k}
\renewcommand\le{\leqslant}
\newcommand\lift{\sigma}
\newcommand{\Linv}{N}
\newcommand{\Low}{L}
\newcommand\mm{m}
\newcommand\MM{M}
\newcommand\nn{n}

\newcommand\NNNN{\mathbb{N}}
\newcommand\pdots{\hspace{0.2ex}{\cdot}{\cdot}{\cdot}\hspace{0.2ex}}
\newcommand\pp{p}
\newcommand\Prod[3]{[#1, #2]_{#3}}

\newcommand\rr{r}
\newcommand\sig[1]{\lift_{\hspace{-0.2ex}#1}^{\null}}
\newcommand{\smallR}{\Sigma}
\renewcommand\ss{s}
\renewcommand\SS{S}
\renewcommand\tt{t}

\newcommand\uu{u}
\def\VR(#1,#2){\vrule width0pt height#1mm depth#2mm}
\newcommand\vv{v}
\newcommand\wAt{\widetilde{\mathrm{A}}_2}
\newcommand\wdots{, ...\HS{0.2},}
\newcommand\ww{w}
\newcommand\WW{W}
\newcommand\WWb{\underline{W}}
\newcommand\XX{X}

\title[Garside families and low elements in Coxeter groups]
{Garside families in Artin--Tits monoids and low elements in Coxeter groups} 

\author{Patrick DEHORNOY}
\address{Laboratoire de Math\'ematiques Nicolas Oresme,
CNRS UMR 6139, Universit\'e de Caen, 14032 Caen, France}
\email{patrick.dehornoy@unicaen.fr}
\urladdr{www.math.unicaen.fr/\!\hbox{$\sim$}dehornoy}

\author{Matthew DYER}
\address{Department of Mathematics, University of Notre Dame, 255 Hurley, Notre Dame, IN 46556, USA}
\email{Matthew.J.Dyer.1@nd.edu}
\urladdr{//math.nd.edu/people/faculty/matthew-j-dyer/}

\author{Christophe HOHLWEG}
\address{Universit\'e du Qu\'ebec \`a Montr\'eal, LaCIM et D\'epartement de math\'ematiques, CP 8888 Succ. Centre-ville, Montr\'eal, Qu\'ebec, H3C 3P8, Canada}
\email{hohlweg.christophe@uqam.ca}
\urladdr{//hohlweg.math.uqam.ca}

\keywords{Coxeter group, Artin--Tits monoid, Garside family, small root, low element, large type, right-angled type, affine type}

\subjclass{20F36, 20F55}

\begin{document}

\begin{abstract}
We show that every finitely generated Artin--Tits group admits a finite Garside family, by introducing the notion of a {\em low element} in a Coxeter group and proving that the family of all low elements in a Coxeter system~$(\WW, \SS)$ with~$\SS$ finite includes~$\SS$ and is finite and closed under suffix and join with respect to the right weak order.
\end{abstract}
 
\maketitle

\selectlanguage{english}
\section{Introduction}\label{S:Intro}

Artin--Tits groups, also known as Artin groups, are those groups defined by relations of the form
\begin{equation}\label{E:Relation}
\ss \tt \ss ... = \tt \ss \tt ...
\end{equation}
where both terms consist of two alternating letters and have the same length. First investigated by J.\,Tits in the late 1960s \cite{Bri}, and then in~\cite{BrS} and~\cite{Dlg}, these groups remain incompletely understood, with many open questions, including the decidability of the Word Problem in the general case~\cite{Cha}. The only well understood case is the one of \emph{spherical type}, which is the case when the associated Coxeter group, obtained by adding the relations $\ss^2 = 1$ to the presentation, is finite. Then a large part of the known results in this case is included in the fact that an Artin--Tits group of spherical type is a Garside group, and the corresponding monoid is a Garside monoid~\cite{Dfx,Dgk}.

At the heart of the properties of an Artin--Tits monoid of spherical type---and more generally of a Garside monoid---lies the fact that every element of the latter admits a distinguished decomposition (``greedy normal form'') involving the divisors of a certain element~$\Delta$ (``Garside element''), in which each entry is in a sense maximal~\cite[Chapter~9]{Eps}. It was recently realized that such distinguished decompositions exist in the more general framework of what was called \emph{Garside families}: whenever $\Gar$ is a Garside family in a left-cancellative monoid~$\MM$ (or category), the mechanism of the greedy normal form works and provides distinguished decompositions with nice properties~\cite{Dif,Garside}. The case of a Garside monoid corresponds to a Garside family consisting of the divisors of a single element~$\Delta$ (``bounded Garside family''), but various examples of unbounded Garside families are now known.

If $\MM$ is an Artin--Tits monoid of non-spherical type, that is, the associated Coxeter group~$\WW$ is infinite, it is well known that $\MM$ is not a Garside monoid: the projection of a possible Garside element to~$\WW$ should be a longest element of~$\WW$, which cannot exist in this case. This however says nothing about possible unbounded Garside families in~$\MM$. In view of effectivity results, the interesting Garside families are the finite ones.
The aim of this note is to announce a proof of the following, previously conjectured, statement, which was supported by partial results and computer experiments.

\begin{thrm}\label{T:ArtinMain}
Every finitely generated Artin--Tits monoid admits a finite Garside family.
\end{thrm}

The proof of Theorem~\ref{T:ArtinMain} relies on translating the problem into the language of Coxeter groups and introducing the new notion of a {\em low element} by looking at the action on the associated root system. Then Theorem~\ref{T:ArtinMain} will follow from the next result, which is of independent interest and seems rich in potential further applications:

\begin{thrm}\label{T:CoxMain}
For every Coxeter system~$(\WW, \SS)$ with~$\SS$ finite, the set of all low elements of~$\WW$ includes~$\SS$ and is finite and closed under join (taken in the right weak order) and suffix.
\end{thrm}

\section{The Artin--Tits problem}\label{S:Artin}

If $\MM$ is a (left-cancellative) monoid and $\ff, \gg$ lie in~$\MM$, one says that $\ff$ \emph{left-divides}~$\gg$ or, equivalently, that $\gg$ is a \emph{right-multiple} of~$\ff$, written $\ff \dive \gg$, if $\ff\gg' = \gg$ holds for some~$\gg'$ in~$\MM$. If there is no nontrivial invertible element in~$\MM$, that is, if $\ff\gg = 1$ holds only for $\ff = \gg = 1$, the left-divisibility relation is a partial ordering that is compatible with multiplication on the left.

\begin{defi}\cite{Dif,Garside}\label{D:Normal}
\rm If $\MM$ is a left-cancellative monoid with no nontrivial invertible element, a \emph{Garside family} of~$\MM$ is a family~$\Gar$ containing~$1$ and such that every element of~$\MM$ admits a $\Gar$-normal decomposition, meaning a finite sequence $(\ss_1 \wdots \ss_\pp)$ satisfying $\ss_1 \pdots \ss_\pp = \gg$ and such that all entries lie in~$\Gar$ and $\forall\ss{\in}\Gar \ \forall\ff{\in}\MM \, (\ss \dive \ff\ss_\ii\ss_{\ii+1} \Rightarrow \ss \dive \ff\ss_\ii)$ holds for every~$\ii < \pp$.
\end{defi}

(Demanding that a Garside family contains~$1$ is not necessary, but it is convenient here, and harmless.) 
The intuition is that a sequence~$(\ss_1 \wdots \ss_\pp)$ is $\Gar$-normal if every entry~$\ss_\ii$ lies in~$\Gar$ and contains as much as possible of the remainder as it can: $\ss_\ii$ is a $\dive$-greatest left-divisor of~$\ss_\ii \pdots \ss_\pp$ lying in~$\Gar$ (whence the word ``greedy'' often used in this context). One shows that, in the context of Definition~\ref{D:Normal}, the $\Gar$-normal decomposition is unique (up to adding or deleting final $1$s) when it exists and that, if $\Gar$ is a Garside family, $(\ss_1 \wdots \ss_\pp)$ is $\Gar$-normal if and only if the simplified condition $\forall\ss{\in}\Gar \, (\ss \dive \ss_\ii\ss_{\ii+1} \Rightarrow \ss \dive \ss_\ii)$ holds for every~$\ii < \pp$.

Various practical characterizations of Garside families are known, depending in particular on the specific properties of the considered monoid. In the case of the Artin--Tits monoid~$\MM$ associated with a Coxeter system~$(\WW, \SS)$, the presentation of~$\MM$ by relations in which both members are words with equal length implies that $\MM$ is strongly Noetherian, meaning that there exits a map~$\lambda: \MM \to \NNNN$ satisfying $\lambda(\gg) \not= 0$ for~$\gg \not=1$ and $\lambda(\gg\hh) \ge \lambda(\gg) + \lambda(\hh)$ for all~$\gg, \hh$. On the other hand, by~\cite[Verk\"urzungslemma]{BrS}, $\MM$ admits conditional right-lcms, that is, any two elements of~$\MM$ that admit a common right-multiple admit a right-lcm (least common right-multiple).

Now, by~\cite[Proposition~3.27]{Dif}, in any left-cancellative monoid~$\MM$ with no nontrivial invertible element that is strongly Noetherian and admits conditional right-lcms, a subfamily~$\Gar$ of~$\MM$ is a Garside family if and only if it contains all atoms of~$\MM$ and is closed under right-lcm and right-divisor. We recall that an element~$\gg$ is called an \emph{atom} if its only left-divisors are~$1$ and~$\gg$, and a family~$\Gar$ is called closed under right-lcm if the right-lcm of two elements of~$\Gar$ belongs to~$\Gar$ when it exists; \ similarly, $\Gar$ is closed under right-divisor if every right-divisor of an element of~$\Gar$ belongs to~$\Gar$, where a \emph{right-divisor} of~$\gg$ is any element~$\ff$ such that $\gg = \gg' \ff$ holds for some~$\gg'$. A direct consequence is that, under the above assumptions, there exists a smallest Garside family in~$\MM$, namely the closure of the atoms under right-lcm and right-divisor \cite[Corollary~3.28]{Dif}. 

Applying this in the case of an Artin--Tits monoid, we obtain:

\begin{Eprop}\label{P:Smallest}
A subfamily~$\Gar$ of an Artin--Tits monoid~$\MM$ is a Garside family if and only if it contains all atoms of~$\MM$ and is closed under right-lcm and right-divisor. In particular, $\MM$ admits a unique smallest Garside family, namely the closure of the atoms under right-lcm and right-divisor.
\end{Eprop}

\begin{coro}\label{C:Main}
An Artin-Tits monoid $\MM$ with atom set~$\SS$ admits a finite Garside family if and only if the closure of~$\SS$ under right-lcm and right-divisor is finite.
\end{coro}

Table~\ref{T:Size} and Proposition~\ref{P:SpecCasesv2} below give some information about the smallest Garside family in a few Artin-Tits monoids. See Figure~\ref{F:AffBraids} for an example.

\begin{figure}[htb]
\centering
\begin{picture}(73,45)(-13,0)
\psset{xunit=0.78mm}
\psset{yunit=0.72mm}
\psset{fillstyle=solid,fillcolor=black}
\cnode(20,30){0.5}{e}
\cnode(30,30){0.5}{a}
\cnode(15,40){0.5}{b}
\cnode(15,20){0.5}{c}
\cnode(20,50){0.5}{ba}
\cnode(35,40){0.5}{ab}
\cnode(5,20){0.5}{cb}
\cnode(5,40){0.5}{bc}
\cnode(20,10){0.5}{ca}
\cnode(35,20){0.5}{ac}
\cnode(30,50){0.5}{aba}
\cnode(0,30){0.5}{bcb}
\cnode(30,10){0.5}{cac}
\cnode[fillcolor=white,linecolor=gray](45,40){0.5}{abc}
\cnode[fillcolor=white,linecolor=gray](45,20){0.5}{acb}
\cnode[fillcolor=white,linecolor=gray](0,50){0.5}{bca}
\cnode[fillcolor=white,linecolor=gray](15,60){0.5}{bac}
\cnode[fillcolor=white,linecolor=gray](15,0){0.5}{cab}
\cnode[fillcolor=white,linecolor=gray](0,10){0.5}{cba}
\cnode(5,0){0.5}{caba}
\cnode(50,30){0.5}{abcb}
\cnode(5,60){0.5}{bcac}
\psset{nodesep=1mm}
\AArrow(e,a)
\BArrow(e,b)
\CArrow(e,c)
\BArrow(a,ab)
\CArrow(a,ac)
\AArrow(b,ba)
\CArrow(b,bc)
\AArrow(c,ca)
\BArrow(c,cb)
\AArrow(ab,aba)
\CArrow(ab,abc)
\BArrow(ba,aba)
\CArrow(ba,bac)
\AArrow(bc,bca)
\BArrow(bc,bcb)
\AArrow(cb,cba)
\CArrow(cb,bcb)
\AArrow(ac,cac)
\BArrow(ac,acb)
\CArrow(ca,cac)
\BArrow(ca,cab)
\BArrow(abc,abcb)
\CArrow(acb,abcb)
\CArrow(bca,bcac)
\AArrow(bac,bcac)
\AArrow(cab,caba)
\BArrow(cba,caba)
\rput[r](e){$1$\ \ }
\rput[l](a){\ \ $\sig1$}
\rput[l](b){\ \ $\sig2$}
\rput[l](c){\ \ $\sig3$}
\rput[r](ab){$\sig1\sig2$\ \ }
\rput[r](ba){$\sig2\sig1$\ \ }
\rput[r](bc){$\sig2\sig3$\ \ }
\rput[r](cb){$\sig3\sig2$\ \ }
\rput[r](ac){$\sig1\sig3$\ \ }
\rput[r](ca){$\sig3\sig1$\ \ }
\rput[r](bcb){$\sig2\sig3\sig2\ \ $}
\rput[l](bac){\ \ $\color{gray}\scriptstyle(\sig2\sig1\sig3)$\ }
\rput[r](bca){$\color{gray}\scriptstyle(\sig2\sig3\sig1)$\ \ }
\rput[l](cac){\ \ $\sig3\sig1\sig3$}
\rput[l](cab){\ \ $\color{gray}\scriptstyle(\sig3\sig1\sig2)$}
\rput[r](cba){$\color{gray}\scriptstyle(\sig3\sig2\sig1)$\ \ }
\rput[l](aba){\ \ $\sig1\sig2\sig1$}
\rput[l](abc){\ \ $\color{gray}\scriptstyle(\sig1\sig2\sig3)$}
\rput[l](acb){\ \ $\color{gray}\scriptstyle(\sig1\sig3\sig2)$}
\rput[l](abcb){\ \ $\sig1\sig2\sig3\sig2$}
\rput[r](bcac){$\sig2\sig3\sig1\sig3$\ \ }
\rput[r](caba){\ $\sig3\sig1\sig2\sig1$\ \ }
\end{picture}
\caption[]{\sf The Cayley graph of the smallest Garside family~$\Gar$ in an Artin--Tits monoid of type~$\wAt$, \ie, the monoid with three generators~$\sig1, \sig2, \sig3$ subject to the relations $\sig1\sig2\sig1 = \sig2\sig1\sig2$, $\sig2\sig3\sig2 = \sig3\sig2\sig3$, $\sig3\sig1\sig3 = \sig1\sig3\sig1$; the Garside family consists of the 16 right-divisors of $\sig1\sig2\sig3\sig2$, $\sig2\sig3\sig1\sig3$, and $\sig3\sig1\sig2\sig1$; the 6 white dots with grey labels do not belong to~$\Gar$, witnessing that $\Gar$ is not closed under left-divisor: $\sig1\sig2\sig3\sig2$ lies in~$\Gar$, whereas its left-divisors $\sig1\sig2\sig3$ and $\sig1\sig3\sig2$ do not.}
\label{F:AffBraids}
\end{figure}

\begin{table}[htb]\centering
\begin{tabular}{c|c|c|c|c|c|c|c|c}
\VR(4,2)\Quad type of $(\WW, \SS)$\Quad\null
&\ spherical\ \null
&\parbox{12mm}{\hfill\VR(2,0)\smash{large}\hfill\null\par\vspace{-4.2mm} \null\hfill\smash{no $\infty$}\VR(0,1)\hfill\null}
&\Quad$\wAt$\Quad\null
&$\Quad\widetilde{\mathrm{A}}_3\Quad\null$
&$\Quad\widetilde{\mathrm{A}}_4\Quad\null$
&$\Quad\widetilde{\mathrm{B}}_3\Quad\null$
&$\Quad\widetilde{\mathrm{C}}_2\Quad\null$
&$\Quad\widetilde{\mathrm{C}}_3\Quad\null$\\
\hline
\VR(4,2)$\card\Garext$&$1$&\VR(0,2.5)$3{\card\SS\choose3}$&$3$&$10$&$35$&$14$&$3$&$12$\\
\hline
\VR(4,2)$\card\Gar$&$\card\WW$&$O((\card\SS)^3)$&$16$&$125$&$1{,}296$&$315$&$24$&$317$
\end{tabular}
\vspace{2mm}
\caption[]{\sf The smallest Garside family~$\Gar$ in the Artin--Tits monoid associated with the Coxeter system~$(\WW, \SS)$; when $\Gar$ is finite, it must consist of all right-divisors of the elements of some minimal finite set~$\Garext$ (the ``extremal elements'').}
\label{T:Size}
\end{table}

\section{Translation of the problem to Coxeter groups}\label{S:Coxeter}

The above considerations admit simple counterparts involving Coxeter groups, which we now explain. Assume that $\MM$ is an Artin--Tits monoid with atom set~$\SS$. Then the quotient~$\WW$ of~$\MM$ obtained by adding the relations $\ss^2 = 1$ to those of~(\ref{E:Relation}) is a Coxeter group. The canonical projection~$\pi$ from~$\MM$ to~$\WW$ is injective on~$\SS$ and (at the expense of identifying $\SS$ with its image under~$\pi$) the pair~$(\WW, \SS)$ is a Coxeter system. By Matsumoto's lemma, mapping a reduced decomposition of an element of~$\WW$ to the element of~$\MM$ admitting that decomposition provides a well defined set-theoretic section~$\lift$ of~$\pi$ from~$\WW$ to~$\MM$, and its image~$\WWb$ is a copy of~$\WW$ inside~$\MM$. 

For~$\ww$ in~$\WW$, we denote by~$\ell(\ww)$ the $\SS$-length of~$\ww$ in~$\WW$, that is, the length of a reduced word for~$\ww$ in $\SS$ (the {\em simple reflections}). Then the product of two elements~$\ff, \gg$ of~$\WWb$ lies in~$\WWb$ if and only if the equality $\ell(\pi(\ff)) + \ell(\pi(\gg)) = \ell(\pi(\ff\gg))$ holds in~$\WW$. We recall that the {\em (right) weak order $\le$} on $\WW$ is defined as follows: let $u,\ww\in W$, then $\uu\le \ww$ holds if and only if a reduced word for $\uu$ is a prefix of a reduced word for $\ww$, if and only if there exists $\vv$ in~$\WW$ satisfying $\ww=\uu\vv$ and $\ell(\ww)=\ell(\uu)+\ell(\vv)$ , see~\cite[Chapter~3]{BjBr}. Now, $(\WW, \le)$ is a complete meet-semilattice~\cite[Theorem 3.2.1]{BjBr}, implying that, if two elements of~$\uu,\vv$ of~ $\WW$ admit a common upper bound with respect to~$\le$, they admit a smallest one called {\em the join}~$\uu\vee\vv$. 

\begin{lemm}\label{L:Dictionary}
Assume that $(\WW, \SS)$ is a Coxeter system and $\MM$ is the associated Artin--Tits monoid.

\ITEM1 The copy~$\WWb$ of~$\WW$ inside~$\MM$ is a Garside family of~$\MM$.

\ITEM2 If $\ff, \gg$ lie in~$\WWb$, then $\ff$ left-divides $\gg$ in~$\MM$ if and only if $\pi(\ff) \le \pi(\gg)$ holds in~$\WW$. Similarly, $\ff$ right-divides~$\gg$ if and only if a reduced word for $\pi(\ff)$ is a \emph{suffix} of a reduced word for~$\pi(\gg)$. 

\ITEM3 If $\ff, \gg$ lie in~$\WWb$, then $\ff$ and $\gg$ of~$\WWb$ have a right-lcm in $\MM$ if and only if $\pi(\ff)\vee\pi(\gg)$ exists in $W$. In this case the right-lcm of $\ff$ and $\gg$ lies in $\WWb$ and is the image under $\lift$ of $\pi(\ff)\vee\pi(\gg)$.
\end{lemm} 

\begin{Proof} 
Point~\ITEM1 follows from~\cite[Proposition~6.27]{Dif}, which says that $\WW$ embeds in the monoid~$\MM'$ generated by~$\WW$ with the relations $\ff\gg = \hh$ for $\ff, \gg, \hh$ satisfying $\ell(\ff) + \ell(\gg) = \ell(\hh)$ and that its image is a Garside family in~$\MM'$. By~\cite{Mic}, the monoid~$\MM'$ is~$\MM$, and the image of~$\WW$ is~$\WWb$.
Next, translating the definition of the left- and right-divisibility relations in~$\MM$ directly gives \ITEM2. 
 Finally, the characterization of a Garside family in an Artin--Tits monoid and~\ITEM1 imply that $\WWb$ is closed under right-lcm in~$\MM$. So, if two elements~$\ff, \gg$ of~$\WWb$ admit a common right-multiple, hence a right-lcm, in~$\MM$, the latter lies in~$\WWb$, and, by~\ITEM2, its projection under~$\pi$ must be the join of~$\pi(\ff)$ and~$\pi(\gg)$, which therefore exists. Conversely, by~\ITEM2, if the join exists, its image under $\lift$ must be the right-lcm of~$\ff$ and~$\gg$ in~$\WWb$. So \ITEM3 is true. 
\end{Proof}

Using the dictionary of Lemma~\ref{L:Dictionary}, we deduce:

\begin{Eprop}\label{P:Translation}
 If $(\WW, \SS)$ is a Coxeter system and $\MM$ is the associated Artin--Tits monoid, the projection of the smallest Garside family of~$\MM$ to~$\WW$ is the smallest subfamily of~$\WW$ that includes~$\SS$ and is closed under join (least common upper bound with respect to the weak order) and suffix.
\end{Eprop}

Thus, in order to prove Theorem~\ref{T:ArtinMain}, it is now sufficient to show that, if $(\WW, \SS)$ is a Coxeter system with~$\SS$ finite, then there exists a finite subset of~$\WW$ that includes~$\SS$ and is closed under join and suffix. 

\section{Low elements in a Coxeter group}\label{S:Low}

The above result will be established by introducing the notion of a {\em low element} in a Coxeter group and showing that the family of all low elements has the expected properties (Theorem~\ref{T:CoxMain}).

From now on, $(\WW, \SS)$ is a fixed Coxeter system with $\SS$ finite. Let $(\Phi,\Delta)$ be a based root system in $(V,B)$ with associated Coxeter system $(W,S)$ as in \cite{HoLaRi14,DyHoRi13}. So, $V$ is a real vector space, $B$ is a symmetric bilinear form on~$V$, and $\Delta$ is a subset of $V$ consisting of one element~$\alpha_\ss$ for each~$\ss$ in~$\SS$ (the \emph{simple} roots). The map sending each $\ss$ in~$\SS$ to the \emph{$B$-reflection} in $\alpha_{\ss}$ extends to a faithful representation of $W$ on~$V$ as the subgroup of the orthogonal group~$O_B(V)$ generated by these reflections. We set $\Phi=W(\Delta)$ (the \emph{roots}), $\Phi^+=\cone(\Delta)\cap \Phi$ (the \emph{positive} roots), and $\Phi^-=-\Phi^+$ (the \emph{negative} roots). Here $\cone(X)$ means the set of all nonnegative linear combinations of elements of~$X$ (the \emph{conic closure} of~$X$). 
 
For~$w$ in~$W$, the {\em (left) inversion set~$\Linv(w)$ of $w$} is $\Phi^+\cap \ww(\Phi^-)$, which is also $\{\alpha\in \Phi^{+}\, | \, \ell(\ss_{\alpha}\ww)<\ell(\ww)\}$. Its cardinality is $\ell(\ww)$. The following properties can be found in, or deduced from,~\cite[Chapter~3]{BjBr} or~\cite{Dy11}.
 
\begin{lemm}\label{L:WeakRoot} 
\ITEM1 For $w$ in~$W$ and $s$ in~$S$ satisfying $\ell(sw)<\ell(w)$, the element $sw$ is a suffix of a reduced word for $w$ and we have $\Linv(w)=\{\alpha_s\}\sqcup s(\Linv( sw))$ (disjoint union).

\ITEM2 The map $\Linv$ is a poset monomorphism from~$(\WW,\le)$ to~$(\mathcal P(\Phi^+),\subseteq )$, and $u\le g$ is equivalent to $\Linv(u)\subseteq \Linv(g)$. 

\ITEM3 For $\uu,\vv$ in~$\WW$ such that $\uu\vee \vv$ exists, $\Linv(u\vee v)=\cone(\Linv(u)\cup \Linv(v))\cap\nobreak \Phi$ holds.
\end{lemm}

In order to show that the language of reduced words of $(W,S)$ is regular, Brink and Howlett introduced in \cite{BrHo93} the notion of dominance order and small roots. The \emph{dominance order} is the partial order~$\domby$ on $\Phi$ such that $\alpha \domby \beta$ holds (``$\beta$ dominates~$\alpha$'') if and only if we have 
$$\forall w\in W\, (w(\beta)\in \Phi^-\Rightarrow w(\alpha)\in\Phi^-).$$
A positive root $\beta$ is called \emph{small}\footnote{These roots are also called \emph{humble} or \emph{elementary} in the literature. We adopt here the terminology of~\cite{BjBr}. See \cite[Notes, p.130]{BjBr} for more details.} when $\beta$ dominates no other positive root than itself, \ie, if we have $\forall \alpha\in\Phi^+\ (\alpha \domby \beta \Rightarrow \alpha=\beta)$. We denote by $\smallR$ the set of small roots. The small roots are characterized recursively by the following lemma. 

\begin{lemm}\cite{BrHo93,BjBr}\label{L:smallchar}
\ITEM1 The set $\Delta$ is included in $\smallR$.

\ITEM2 For every~$\beta$ in $\Phi^{+}\setminus \Delta$, there exists $\alpha$ in~$\Delta$ satisfying $\ell(\ss_{\alpha}\ss_{\beta}\ss_{\alpha})=\ell(\ss_{\alpha})-2$, or equivalently, $B(\alpha,\beta)>0$. Then, for every such $\alpha$, one has $\beta\in \smallR$ if and only if $s_{\alpha}(\beta)$ lies in~$\smallR$ and $B(\alpha,\beta)<1$ holds.
\end{lemm}

\begin{thrm}[Brink-Howlett~\cite{BrHo93}]\label{thm:BrHo} The set $\smallR$ is finite. 
\end{thrm}

The following result is a restatement of a special case of Propositions 1.4 and 3.6 in~\cite{Dy94}:

\begin{lemm}\label{L:extremeray} 
For $\ww$ in~$\WW$, let $\Linv^{1}(\ww):=\{\alpha\in \Phi^{+}\mid \ell(\ss_{\alpha}\ww)=\ell(\ww)-1\} \subseteq \Linv(\ww)$. Then $\Linv^{1}(\ww)$ is the set of all $\alpha$ in~$\Phi^{+}$ such that the cone of~$\{\alpha\}$ is an extreme ray of the polyhedral cone spanned by~$\Linv(\ww)$.
\end{lemm}

\begin{defi} 
\rm An element $\ww$ of~ $\WW$ is {\em low} if $\Linv^{1}(\ww) \subseteq \smallR$ holds, \ie, if we have $\Linv(\ww) = \cone(A) \cap \Phi$ for some some family of small roots~$A$. We denote by $\Low$ the set of low elements of~$\WW$. 
\end{defi}

We can now sketch the proof of Theorem~\ref{T:CoxMain}; a complete proof of this theorem can be found in~\cite{DyHo}.

\begin{sketch} 
First, the fact that $\Low$ is finite follows from Theorem~\ref{thm:BrHo} and Lemma~\ref{L:WeakRoot}\ITEM2: there is only a finite number of subsets of~$\smallR$, hence a finite number of low elements, since the map $\Linv$ is injective. Then, the fact that $\Low$ includes~$\SS$ follows from the fact that, for~$\ss$ in~$\SS$, we have $\Linv(\ss)=\{\alpha_\ss\}\subseteq \smallR$. Now, assume that we have $\Linv(\uu)=\cone(A)\cap \Phi$ and $\Linv(\vv)=\cone(B)\cap \Phi$ with $A,B \subseteq \smallR$. By definition of the conic closure, we have $\cone(\cone(A)\cup\cone(B))=\cone(A\cup B)$. By Lemma~\ref{L:WeakRoot}\ITEM3, we deduce 
$$\cone(\Linv(\uu\vee \vv))\cap\Phi=\cone(\cone(A)\cup\cone(B))\cap \Phi=\cone(A\cup B)\cap \Phi:$$
 as $A\cup B$ is included in~$\smallR$, we conclude that $\uu \vee \vv$ lies in~$\Low$, so $\Low$ is closed under join. 

The difficult part is to show that $\Low$ is closed under suffix, and here we only give a sketch of the proof. Recall first that a {\em maximal rank~$2$ root subsystem of $\Phi$} is a set~$\Phi'$ of the form $\Phi'=P\cap \Phi$ where $P$ is a plane in $V$ intersecting $\Phi^+$ in at least two roots. The cone spanned by~$\Phi' \cap\Phi^+$ has then a basis $\Delta'$ of cardinality $2$ included in~$\Phi'\cap \Phi^+$, and then one has $P\cap \Phi^{+}=\cone(\Delta')\cap \Phi$. 

\ITEM1 We start by showing that $\smallR$ is {\em bipodal}, meaning that, for every small root $\beta$ and for every maximal rank~$2$ root subsystem $\Phi'$ of $\Phi$ with basis $\Delta'$ satisfying $\beta\in\Phi'\setminus\Delta'$, we must have $\Delta'\subseteq \smallR$. To prove this, we note that, for every~$\alpha$ in ~$\Delta$ satisfying $B(\alpha,\beta)>0$, the reflection subgroup generated by reflections in $\Delta'\cup\{\alpha\}$ is of rank at most three and $\beta$ is a small root for its corresponding root subsystem. Using this observation and Lemma \ref{L:smallchar}, one reduces by induction on $\ell(\ss_{\beta})$ to the case of root systems of rank three. Then one checks the result in rank three using the explicit descriptions of small roots in~\cite{BrDom}.

\ITEM2 Now, consider $\ww$ and~$\ss$ as in Lemma~\ref{L:WeakRoot}\ITEM1. Write $\ss=\ss_{\alpha}$ with $\alpha\in \Delta$. For every $\beta$ in~$\Phi^{+}\setminus \{\alpha\}$, let $\ff_{\alpha}(\beta)$ be the simple root different from $\alpha$ in the standard simple system of the maximal rank two root subsystem~$\Phi\cap P$, where $P$ is the plane spanned by $\alpha$ and $\beta$: in other words, we have $\ff_{\alpha}(\beta)\in \Phi^{+}$ and $P\cap \Phi^{+}=\cone(\{\alpha, \ff_{\alpha}(\beta)\})\cap \Phi^{+}$. Then we show that $\Linv^{1}(\ss\ww)$ is included in~$\{ \ss(\beta), \ff_{\alpha}(\beta) \mid \beta\in \Linv^{1}_{\ww}\setminus\{\alpha\}\}$. To prove this, one uses Lemma~\ref{L:extremeray} to reformulate it as a statement in terms of Bruhat order, and checks that statement using standard properties from \cite{Dy90} of cosets of (maximal dihedral) reflection subgroups. 

\ITEM3 Finally, to show that $\Low$ is closed under suffix, it is enough to show that, for~$\ww$ in~$\Low$ (\ie, for $\Linv^{1}(\ww) \subseteq \smallR$) and $\ss$ in~$\SS$ satisfying $\ell(\ss\ww)<\ell(\ww)$, the element~$\ss\ww$ also lies in~$\Low$. Write $\ss=\ss_{\alpha}$ with~$\alpha$ in~$\Delta$. By~\ITEM2, it is sufficient to show that $\ss(\beta)$ and $\ff_{\alpha}(\beta)$ are small for every $ \beta$ in~$\Linv^{1}(\ww)\setminus\{\alpha\}$. But, by~\ITEM1, $\ff_{\alpha}(\beta)$ is small since $\beta$ is. On the other hand, assume that $\ss(\beta)$ is not small. Then, Lemma \ref{L:smallchar} implies $B(\alpha,\beta)\le -1$ . But then the subgroup generated by $\ss$ and $\ss_{\beta}$ is infinite dihedral with $\alpha$ and $\beta$ as its simple roots, so its positive system $\Phi^{+}\cap \cone(\{\alpha,\beta\})$, which is infinite, must be included in $\Linv(\ww)$, which is finite. This contradiction shows that $\ss(\beta)$ must be small, and completes the proof. 
\end{sketch}

In \cite{DyHo}, we introduce and study a more general notion, the {\em $\nn$-low elements}, that are defined using a notion of $\nn$-small roots. The $0$-low elements are the low elements as defined in this text. We conjecture that, for every~$\nn$, the $\nn$-low elements give rise to a Garside family in the associated Artin--Tits monoid. We know that they form a finite set, closed under join. To conclude, it would suffice to prove that the set of $\nn$-small roots is bipodal for every $\nn$, which we conjecture in general and prove in some cases including affine Weyl groups.

\section{Descriptions of $\pi(\Gar)$ and $\Low$ in some special cases}
 
We keep the same notation, and describe the image~$\pi(\Gar)$ of the smallest Garside family~$\Gar$ of~$\MM$ in a few cases. By Proposition~\ref{P:Translation}, $\pi(\Gar)$ is the closure of~$\SS$ under join and suffix in~$\WW$. We denote the Coxeter matrix of $(\WW,\SS)$ as $(\mm_{\ss,\tt})_{\ss,\tt\in \SS}$, and write $\Prod\ss\tt\kk$ for the alternating product $\ss\tt\ss\tt\pdots$ with $\kk$ factors, $\kk\ge 1$. It is well-known that the {\em standard dihedral parabolic} subgroup $W_{\{\ss,\tt\}}$ consists of the identity together with the elements $\Prod\ss\tt\kk$ and $\Prod\tt\ss\kk$, $k\ge 1$. Moreover $W_{\{\ss,\tt\}}$ is finite if and only if $\mm_{\ss,\tt}$ is finite and in this case the longest element is $\Prod\tt\ss{m_{\ss, \tt}}=\Prod\ss\tt{m_{\ss, \tt}}$.

\begin{Eprop}\label{P:SpecCasesv2}
\ITEM1 If $\MM$ is an Artin--Tits monoid of spherical type, then we have $\pi(\Gar) =\Low=\WW$.

\ITEM2 If $\MM$ is an Artin--Tits monoid of large type (\ie, $\mm_{\ss, \tt} \ge3$ holds for all~$\ss\neq \tt$), then we have $\pi(\Gar) =\Low=\XX$, where $\XX$ is the union of all finite standard parabolic subgroups of $\WW$ (each being of rank at most two) together with all elements $t\Prod\rr\ss{m_{r,s}}$ with $r,s,t$ distinct in~$S$ and $m_{r,s}$, $m_{s,t}$, $m_{t,r}$ all finite. 
 
\ITEM3 If $\MM$ is a right-angled Artin--Tits monoid (\ie, $\mm_{\ss, \tt} \in \{2, \infty\}$ holds for all~$\ss \not=\tt$), then we have $\pi(\Gar) = \Low=\XX$, where $\XX$ is the union of all finite standard parabolic subgroups of $W$ (which are of the form $W_{I}$ where $\II$ is a set of pairwise commuting simple reflections). 
\end{Eprop}

\begin{sketch} 
First, $\lift(\Low)$ is a Garside family in~$\MM$ by Theorem~\ref{T:CoxMain}, which implies $\Gar \subseteq \lift(\Low)$, whence $\pi(\Gar) \subseteq \Low$ in every case. Hence, for~\ITEM1, it suffices to show $\WW \subseteq \pi(\Gar)$ and, for~\ITEM2 and~\ITEM3, it suffices to show $\XX\subseteq \pi(\Gar)$ and $\Low\subseteq \XX$. 

\ITEM1 Here, $\pi(\Gar)$ contains the join of all elements of $\SS$, which is the longest element $\ww_{0}$ of $\WW$, and every element of $\WW$ is a suffix of $\ww_{0}$. We deduce $\WW\subseteq \pi(\Gar)$, as required. (Note that, if $W$ is infinite, then any finite standard parabolic subgroup~$W_I$ generated by a subset~$I$ of~$S$ satisfies $W_I\subseteq \pi(\Gar)$ by the same arguments.)

\ITEM2 First, $S\subseteq \pi(\Gar)$ holds by definition. Next, for $\rr, \ss$ distinct in~$\SS$ with $\mm_{\rr, \ss}$ finite, the subgroup $W_{\{\rr,\ss\}}$ is finite and by the remark above we have $ W_{\{\rr,\ss\}} \subseteq\pi(F)$. Finally, for $\rr,\ss,\tt$ pairwise distinct in~$\SS$ with $\mm_{\rr,\ss}$, $\mm_{\ss,\tt}$ and $\mm_{\tt,\rr}$ all finite, as just seen, $\tt\rr$ and $\tt\ss$ lie in~$\pi(\Gar)$, hence so does their join, which is $\tt\Prod\rr\ss{\mm_{\rr, \ss}}$. This shows $\XX\subseteq \pi(\Gar)$.

Now we show $\Low\subseteq \XX$. First, by~\cite{BrDom}, the full subgraph of the Coxeter graph on the support of a small root contains no cycle or infinite bond. Hence, in large type, the small reflections (meaning the reflections in a small root) are precisely the reflections in the finite standard parabolic subgroups. Assume that $\rr,\ss,\tt$ are pairwise distinct in~$\SS$. We claim that an element of $\Low$ cannot admit a reduced expression of the form $\ww=\uu \tt\Prod\rr\ss\kk$ with $\uu \not= \tt \in \SS$ and $2\le \kk\le\mm_{\rr,\ss}$. Indeed, assume $\ww \in \Low$. For any reduced expression $\rr_\nn \pdots \rr_1$ of~$\ww$ and~$1\leq \ii\leq n$ 
we define $\tt_{\ii}:=\rr_{\nn}\pdots \rr_{\ii+1}\rr_{\ii}\rr_{\ii+1}\pdots \rr_\nn$, a reflection with $\ell(\tt_{\ii}\ww)<\ell(\ww)$. By Lemma~\ref{L:extremeray}, $\tt_{\ii}$ is a small reflection if $\ell(\tt_{\ii}\ww)=\ell(\ww)-1$ holds. So, here, $\tt_{\nn-2}$, \ie, $\uu\tt\rr\tt\uu$, must be a small reflection, which forces $\uu=\rr$ (and $\mm_{\rr,\tt}<\infty$). Also, $\tt_{1}$ must be a small reflection. Now, we have $\tt_{1}=\uu\tt\vv\tt\uu=\rr\tt\vv\tt\rr$ with $\vv=\Prod\rr\ss{2\kk-1}$, and $\vv$ is a reflection of~$\WW_{\{\rr,\ss\}}$ unequal to~$\rr$. Since $\rr$, $\ss$ and $\tt$ are pairwise non-commuting, Matsumoto's Lemma implies, first, $\ell(\tt\vv\tt)=\ell(\vv)+2$ and, then, $\ell(\tt_{1})=\ell(\tt\vv\tt)+2$, by considering the cases $\vv=\ss$ and $\ell(\vv)>1$ separately. Hence the smallest standard parabolic subgroup containing~$\tt_{1}$ is $\WW_{\{\rr,\ss,\tt\}}$ of rank $3$, a contradiction.

Similar (but simpler) arguments show that an element of~$\Low$ cannot admit a reduced expression of the form $\tt \Prod\rr\ss\kk$ with $\rr,\ss, \tt$ distinct in~$\SS$ and $2\le \kk < \mm_{\rr,\ss}$, or $\Prod\rr\ss\kk$ with $\rr, \ss$ distinct in~$\SS$ and $2\le k<\nobreak\mm_{\rr,\ss}=\infty$. Now, every element~$w$ of~$W$ has a unique decomposition $w=uv$ with $u\in W$ satisfying $\ell(us)=\ell(ur)>\ell(u)$ and $v\in W_{\{\rr,\ss\}}$, see~\cite[Proposition 2.4.4 (i)]{BjBr}. So, as $\Low$ is closed under suffix, no element of any of the three types excluded above can be a suffix of an element of~$\Low$, and $\Low\subseteq \XX$ easily follows. 

\ITEM3 The proof is similar to (and simpler than) that of~\ITEM2. First, since small roots cannot have any infinite bonds in their supports, the set of small reflections is precisely~$S$. Then, for $\II\subseteq \SS$ consisting of commuting simple reflections, the parabolic subgroup $W_I$ is finite. Hence such $W_{I}$ are contained in~$\pi(\Gar)$, which implies $\XX\subseteq \pi(\Gar)$. 

For $\Low\subseteq \XX$, suppose $\ww=\uu\rr_{1}\pdots \rr_{\nn}$ where $\uu,\rr_{1} \wdots \rr_{\nn}$ are distinct pairwise commuting simple reflections such that $\uu$ does not commute with all of $\rr_{1} \wdots \rr_{\nn}$, and let $\ii$ be minimal with $\uu\rr_{i}\neq \rr_{\ii}\uu$. If $\ww$ were low, arguing as in~\ITEM2, we would deduce that $\uu\rr_{1}\pdots \rr_{i}\pdots \rr_{1}\uu$, \ie, $\uu\rr_{i}\uu$, is a small reflection, so lies in $S$, leading to a contradiction with $\uu\rr_{i}\neq \rr_{\ii}\uu$. Hence $\uu\rr_{1}\pdots \rr_{\nn}$ is not low and, again, one easily deduces $\Low\subseteq X$.
\end{sketch}

In type $\widetilde{\mathrm{C}}_2$, with Coxeter graph 
 $ \xymatrix{\sigma_{1}\ar@{=}[r]&\sigma_{2}\ar@{=}[r]&\sigma_{3}}$, one finds $\vert \Low\vert =25$ and $\vert \pi(\Gar)\vert = \vert \Gar\vert=24$: here $\sig1\sig3\sig2$ is low, but does not lie in~$\pi(\Gar)$. Hence $\pi(\Gar) = \Low$ need not hold in general. 
 
However, for type~$\widetilde{\mathrm{A}}_\nn$, the equality $\pi(\Gar) = \Low$ holds and we have $\vert \Low \vert = (\nn+ \nobreak2)^\nn$. Indeed, while preparing this manuscript, the third author (CH), together with P.\, Nadeau and N.\, Williams, built two automata recognizing the language of reduced words with respective state sets $\pi(\Gar)$ and $\Low$. It is easy to see that $\pi(\Gar)$ and $\Low$ inject into the state set of the {\em canonical automaton} defined by Brink and Howlett, see~\cite[p.120]{BjBr}. Now, in type~$\widetilde{\mathrm A}_\nn$, Eriksson showed that the latter is minimal \cite[p.125]{BjBr}, so $\pi(\Gar)$ and $\Low$ must share its cardinality, which is~$(\nn+2)^\nn$. A more direct proof would be desirable.

\bibliographystyle{plain}

\end{document}